\documentclass[12pt]{article}
\usepackage{epsfig,graphicx,color}
\usepackage{amsmath}
\usepackage{latexsym}
\usepackage{amstext}
\usepackage[normalem]{ulem}

\newcommand{\mb}[1]{ \mbox{\boldmath$#1$}}
\newcommand{\ds}{\displaystyle}
\newcommand{\beq}{\begin{eqnarray}}
\newcommand{\eeq}{\end{eqnarray}}
\newcommand{\beqq}{\begin{eqnarray*}}
\newcommand{\eeqq}{\end{eqnarray*}}
\newcommand{\p}{\partial}

\newcommand{\eps}{\varepsilon}
\newcommand{\x}{\mbox{\boldmath$x$}}
\newcommand{\h}{\mbox{\boldmath$h$}}

\newcommand{\y}{\mbox{\boldmath$y$}}

\font\bb=msbm10 at 12pt
\def\rR{\hbox{\bb R}}

\def\r{\right}
\def\l{\left}

\begin{document}
\pagestyle{plain}
\begin{center}
{\large \textbf{{Nonlinear Filtering with Optimal MTLL}}
\\[5mm]
E. Fischler
\footnote{Department of Systems, School of Electrical
Engineering, The Iby and Aladar Fleischman Faculty of Engineering,
Tel-Aviv University, Ramat-Aviv Tel-Aviv 69978, Israel.},
Z. Schuss
\footnote{Department of Mathematics, Tel-Aviv University, Tel-Aviv
69978, Israel.}}
\end{center}
\date{}

\begin{abstract}
\noindent We consider the problem of nonlinear filtering of
one-dimensional diffusions from noisy measurements. The filter is
said to lose lock if the estimation error exits a prescribed region.
In the case of phase estimation this region is one period of the
phase measurement function, e.g., $[-\pi,\pi]$. We show that in the
limit of small noise the causal filter that maximizes the mean time
to loose lock is Bellman's minimum noise
energy filter.
\end{abstract}

\section{Introduction}
Optimal filtering theory defines different optimality criteria, such
as minimizing the conditional mean square estimation error (MSEE),
given the measurements \cite{Jazwinsky}, maximizing the a posteriori
probability (MAP) density function (pdf) of the signal, given the
measurements \cite{VanTrees}, Bellman's criterion of minimum noise
energy  (MNE) \cite{Bellman}, \cite{Lee}, \cite{Lee1}, and more. In
problems of phase estimation, that lead to loss of lock and cycle
slips, an important optimality criterion is maximizing the mean time
to lose lock (MTLL) or to exit a given region, which is also a well
known control problem \cite{Fleming}, \cite{ZZ}, \cite{ADS},
\cite{MR}. Approximation methods for finding the various optimal
filters have been devised for problems with small noise, including
large deviations and WKB solutions of Zakai's equation, the extended
Kalman filter (EKF) \cite{Katzur}, \cite{Picard86}\cite{Picard91},
\cite{Hijab} and others. The EKF and WKB approximations produce explicit suboptimal
finite-dimensional filters, which in case of phase estimation are
the well known phase trackers, such as the phase locked loop (PLL),
delay locked loop (DLL), angle tracking loops, and so on
\cite{Stensby}. The MSEE in these phase trackers is asymptotically
optimal \cite{Katzur}, \cite{Hijab}.

The suboptimal phase trackers are known to lose lock (or slip
cycles) \cite{Stensby}. The MTLL in these filters is simply the mean
first passage time (MFPT) of the estimation error to the boundary of
the lock region. The MFPT from an attractor of a dynamical system
driven by small noise has been calculated by large deviations and
singular perturbation methods \cite{FW}, \cite{MS}, \cite{book}, and
in particular, for the PLL \cite{BS}. The MTLL in particle filters
for phase estimation was found in \cite{FB}. It has been found
recently that minimizing the MNE leads to a finite, yet much longer
MTLL than in the above mentioned phase estimators \cite{Doron},
\cite{Ehud}. This raises the question of designing a causal (or
noncausal) phase estimator with maximal MTLL.

The MTLL is the fundamental performance criterion in phase tracking
and synchronization systems. Thus, for example, a phase tracking
system is considered locked, as long as the estimation error
$e(t)=x(t)-\hat x(t)$ is in $(-\pi,\pi)$. When the error exceeds
these limits, the estimation is said to be unlocked, and the system
relocks on an erroneous equilibrium point, with  a deviation of
$2\pi$. Another example is an automatic sight of a cannon. The sight
is said to be locked on target if the positioning error is somewhere
in between certain limits. Similar problems, in which the
maximization of exit time is an optimality criterion, were
considered by several authors \cite{ZZ}, \cite{ADS}, \cite{MR}. In
\cite{ZZ}, a simpler filtering problem is considered, in which the
error $e(t)$ is measured, rather than the state variable $x(t)$. It
is solved under the further assumption of a linear measurement
inside a domain. In \cite{ADS}, \cite{MR} the state process is
controlled through its drift, rendering it a control rather than a
filtering problem.

In this paper we show that for small noise the maximum MTLL filter
is Bellman's MNE filter \cite{Lee}. It follows that the result of
\cite{Ehud} for the MTLL of the optimal MNE phase filter, is
asymptotically an upper bound for any other filtering scheme. In
view of the results of \cite{Ehud}, the potential gain of the
optimal MNE filter over the first order EKF-PLL is $12$ dB.

\section{Formulation}\label{sec:Form}
An important class of filtering problems with small
measurements noise can be reduced to the model of a diffusion
process
 \beq\label{proc_eq}
    dx(t)&=&m(x,t)\,dt+\eps\sigma\,dw(t),
 \eeq
measured in a noisy channel
 \beq\label{meas_eq}
    dy(t)&=&h(x,t)\,dt+\eps\rho\,d\nu(t),
 \eeq
where $m(x,t)$ and $h(x,t)$ are possibly nonlinear, continuous
functions. The processes $w(t)$ and $\nu(t)$ are independent
standard Brownian motions, and $\eps$ is a small parameter. If
$m(x,t)$ is a linear function and the noise in (\ref{proc_eq}) is
not small, an appropriate scaling of time and dependent variables
scales the small measurements noise into the diffusion equation as
well, giving the canonic system (\ref{proc_eq}), (\ref{meas_eq})
\cite{Doron}. The optimal filtering problem is to find a causal
estimator $\hat x(t)$ of $x(t)$, given the measurements
$y_0^t=\{y(s)\,:\, 0\leq s\leq t\}$, such that the mean first time
the error signal,
 \beq
    e(t)=x(t)-\hat x(t),\label{exhat}
 \eeq
leaves a given lock domain $L\subset\rR$, is maximal. More
specifically, for any adapted function $\hat x(t)\in{\cal C}(\rR^+)$
(measurable  with respect to the filtration generated by $y(t)$), we
define an error process by (\ref{exhat}) and the first time to lose
lock by
 \beq
     \tau=\inf\left\{t\,:\,e(t)\in\p L\right\}.\label{tauL}
 \eeq
The optimal filtering problem is to maximize $E[\tau\,|\,y_0^\tau]$
(see definition (\ref{Etauy}) below) with respect to all
adapted continuous functions $\hat x(t)$. For example, if
$h(x,t)=\sin x$ in a phase estimation problem, then $L=(-\pi,\pi)$
and lock is lost when $e(t)=\pm\pi$.

We can rewrite the model equations (\ref{proc_eq}), (\ref{meas_eq})
in terms of the error process $e(t)$ as
 \beq
    de(t)&=&M_{\hat x}(e(t),t)\,dt+\eps\sigma\,dw(t)\label{err_eq}\\
    &&\nonumber\\
    dy(t)&=&H_{\hat x}(e(t),t)\,dt+\eps\rho\,d\nu(t)\label{err_meas}
 \eeq
where
 \beqq M_{\hat x}(e(t),t)&=&m(\hat x(t)+e(t))-\dot{\hat x}(t)\\
 &&\\
     H_{\hat x}(e(t),t)&=&h(\hat x(t)+e(t)),\eeqq
and the filtering problem is to find $\hat x(t)$, such that
$E[\tau\,|\,y_0^\tau]$ is maximal.

The survival probability of a trajectory $(e(t),y(t))$ of
(\ref{err_eq}) with absorption at $\p L$ and (\ref{err_meas}) can be
expressed in terms of the pdf $p(e,y,t\,\big|\,\xi,\eta,s)$ of the
two-dimensional process with an absorbing boundary condition on $\p
L$. It is the solution of the Fokker-Planck equation (FPE)
 \beq
    \frac{\p p(e,y,t\,\big|\,\xi,\eta,s)}{\p t}&=&-\frac{\p
    M_{\hat x}(e,t)p(e,y,t\,\big|\,\xi,\eta,s)}{\p e}-\frac{\p
    H_{\hat x}(e,t)p(e,y,t\,\big|\,\xi,\eta,s)}{\p y}+\nonumber\\
    &&\nonumber\\
    && \frac{\eps^2\sigma^2}{2}\frac{\p^2
    p(e,y,t\,\big|\,\xi,\eta,s)}{\p
    e^2}+\frac{\eps^2\rho^2}{2}\frac{\p^2
    p(e,y,t\,\big|\,\xi,\eta,s)}{\p y^2}\label{2DFPE}
 \eeq
for $e,\xi\in L,\ y,\eta\in\rR$, with the boundary and initial
conditions
 \beq
    p(e,y,t\,|\,\xi,\eta,s)&=&0\qquad\qquad\qquad\quad\mbox{for}\quad e\in\p L,\
    y\in\rR,\ \xi\in L,\ \eta\in\rR\label{BC}\\
    &&\nonumber\\
    p(e,y,s\,|\,\xi,\eta,s)&=&\delta(e-\xi,y-\eta)\ \quad\mbox{for}\quad
    e\in L,\ y\in\rR,\ \xi\in L,\ \eta\in\rR\label{IC}.
 \eeq
The pdf is actually the joint density and probability function
$p(e,y,t\,|\,\xi,\eta,s)=\Pr\{e(t)=e,y(t)=y,\tau>t\,|\,\xi,\eta,s\}$
and thus the survival probability is
 \beq
\Pr\{\tau>t\,|\,\xi,\eta,s\}=S_{e(\cdot),y(\cdot)}(t)&=&
\int_{L}\int_{\rR} p(e,y,t\,|\,\xi,\eta,s)\,de\,dy,\label{Sxy}
 \eeq
and it decays in time.

\section{Simulation with particles}
To simulate the filtering problem on a finite interval $0\leq t\leq
T$, we discretize (\ref{proc_eq}), (\ref{meas_eq}) on a sequence of
grids
$$\left\{t_i=i\Delta t,\quad i=0,1,\ldots,N,\quad\Delta t=\frac{T}{N}\right\},$$
and define discrete trajectories by the Euler scheme
   \beq
    x_N(t_{i+1})&=&x_N(t_i)+\Delta t\,m(x_N(t_i),t_i)+
    \eps\sigma\,\Delta w(t_i)
    \label{dxN}\\ &&\nonumber\\
    y_N(t_{i+1})&=&y_N(t_i)+\Delta t\,h\left(x_N(t_i),t_i\right)
    +\eps\rho\,\Delta\nu(t_i), \label{dyN}
 \eeq
for $i=0,1,\ldots,N-1$, where $\Delta w(t_i)$ and $\Delta
\nu(t_i)$ are independent zero mean Gaussian random variables with
variance $\Delta t$. The discretized version of (\ref{err_eq}),
(\ref{err_meas}) is
 \beq
 e_N(t_{i+1})&=&e_N(t_i)+\Delta t\,M_{\hat x}\left(e_N(t_i),t_i\right)+
 \eps\sigma\,\Delta w(t_i)
    \label{de}\\ &&\nonumber\\
    y_N(t_{i+1})&=&y_N(t_i)+\Delta t\,H_{\hat x}\left(e_N(t_i),t_i\right)
    +\eps\rho\,\Delta\nu(t_i).
    \label{dy}
 \eeq
Given an observed trajectory $\{y_{N}(t_i)\}_{i=0}^N$, we
sample $n$ trajectories
$\left\{\{x_{j,N}(t_i)\}_{i=0}^N\right\}_{j=1}^n$, according to the
scheme (\ref{dxN}), which produce error trajectories
$\left\{\{e_{j,N}(t_i)\}_{i=0}^N\right\}_{j=1}^n$, and determine
their first exit times from $L$, denoted
$\left\{\tau_{j,N}\right\}_{j=1}^n$ (we set $\tau_{j,N}=T$ if
$\{e_{j,N}(t_i)\}_{i=0}^N$ does not exit $L$ by time $T$)
\cite{Crisan2003}, \cite{Amblard}, \cite{Fischler},
\cite{Arulampalam}, \cite{elMoral_Guionnet}. The conditional MTLL is
defined on the ensemble by
 \beq
 &&E\left[\tau_N\wedge T\,|\,\left\{y_N(t_i),\,i=0,1,\ldots,\ds\left[\frac{\tau_{N}}
 {\Delta t}\right]\wedge N\right\}\right]=\label{EtauN}\\
 &&\nonumber\\
 &&\frac{\ds\sum_{j=1}^n(\tau_{j,N}\wedge T)\exp
 \left\{\ds\frac{1}{\varepsilon^2\rho^2}\ds\sum_{k=0}^{\left[\ds\frac{\tau_{j,N}\wedge T}{\Delta t}\right]}\left[H(e_{j,N}(t_{k-1}),t_{k-1})\Delta y_{k,N}-\frac12
 H^2(e_{j,N}(t_{k-1}),t_{k-1})\Delta t\right]\right\}}{\ds\sum_{j=1}^n\exp
 \left\{\ds\frac{1}{\varepsilon^2\rho^2}\ds\sum_{k=0}^{\left[\ds\frac{\tau_{j,N}\wedge T}{\Delta t}\right]}\left[H(e_{j,N}(t_{k-1}),t_{k-1})\Delta y_{k,N}-\frac12
 H^2(e_{j,N}(t_{k-1}),t_{k-1})\Delta t\right]\right\}}.\nonumber
 \eeq
We define
 \beq
 E[\tau\,|\,y_0^{\tau}]=
 \lim_{T\to\infty}\lim_{n\to\infty}\lim_{N\to\infty}
 E\left[\tau_N\wedge T\,|\,\left\{y_N(t_i),\,i=0,1,\ldots,\ds\left[\frac{\tau_{N}}
 {\Delta t}\right]\wedge N\right\}\right].\label{Etauy}
 \eeq
The conditional MTLL $E[\tau\,|\,y_0^{\tau}]$ is a random variable
on the $\sigma$-algebra ${\cal F}_\infty=\ds\bigcup_{t>0}{\cal
F}_t$, where ${\cal F}_t$ is the $\sigma$-algebra generated by the
measurements process $y(\cdot)$ up to time $t$. Our purpose is to
find $\hat x(t)$ that maximizes $E[\tau\,|\,y_0^\tau]$ in the class
of continuous adapted functions.
\section{The joint pdf of the discrete process}
The pdf of a trajectory of $(e_N(t),y_N(t))$ is the Gaussian
 \beq
 &&p_N(e_1,e_2,\ldots,e_N;y_1,y_2,\ldots,y_N;t_1,t_2,\dots,t_N)=
 \prod_{k=1}^N\left[\frac{\ds\exp\l\{-\frac{{\cal
    B}_k(\x_k,\x_{k-1})}{2\varepsilon^2\Delta t }\right\}}
    {2\pi\varepsilon^2\rho\sigma\Delta t }\right],\label{peNyN}
    \eeq
where the exponent is the quadratic form
 \beqq
    {\cal B}_k(\x_k,\x_{k-1})&=&\left[\x_k-\x_{k-1}-\Delta t\mb{a}_{k-1}\right]^T
    \mb{B}\left[\x_k-\x_{k-1}-\Delta t\mb{a}_{k-1}\right],
 \eeqq
such that
 \beqq
     \x_k=\left[\begin{array}{l}e_k\\y_k\end{array}\right],\quad
     \mb{a}_k = \left[\begin{array}{l}M_{\hat x}(e_k,t_k)\\H_{\hat x}(e_{k},t_k)\end{array}\right],\quad\mb{B}=
     \begin{bmatrix} \sigma^{-2} & 0 \\ 0 & \rho^{-2} \\
     \end{bmatrix}.
 \eeqq

The Wiener path integral \cite{Schulman}, \cite{Freidlin},
\cite{Kleinert}, \cite{unidirect}, \cite{ZakaiUstunel}
 \beq
    &&p(e,y,t\,|\,\xi,\eta,s)=\label{pxyt}\\
    &&\nonumber\\
    &&\lim_{N\to\infty} \underbrace{\int_{L}de_1\int_{L}de_2\cdots
    \int_{L}de_{N-1}}_{N-1}\underbrace{\int_{\rR}dy_1\int_{\rR}dy_2\cdots
    \int_{\rR}dy_{N-1}}_{N-1}\times\nonumber \\
    &&\nonumber\\
    &&\prod_{k=1}^N\left[\frac{\ds\exp\l\{-\frac{{\cal B}_k(\x_k,\x_{k-1})}{2\varepsilon^2\Delta t }\right\}}
    {2\pi\varepsilon^2\rho\sigma\Delta t }\right],\nonumber
 \eeq
with $e_N=e,\ y_N=y,\ e_0=\xi,\ y_0=\eta$, is the solution of the
FPE (\ref{2DFPE}) with the boundary and initial conditions
(\ref{BC}) and (\ref{IC}).

The pdf (\ref{peNyN}) can be written as
 \beq
 &&p_N(e_1,e_2,\ldots,e_N;y_1,y_2,\ldots,y_N;t_1,t_2,\dots,t_N)=\label{break}\\
 &&\nonumber\\
 &&\prod_{k=1}^N\left[\frac{1}{{\sqrt{2\pi\Delta t }}\,\varepsilon\sigma}\exp\l\{
 -\frac{[e_k-e_{k-1}-\Delta t M_{\hat x}(e_{k-1},t_{k-1})\,]^2}
    {2\eps^2\sigma^2\Delta t }\r\}\times\right.\nonumber\\
   &&\nonumber\\
 &&\left. \exp\l\{\frac{1}{\eps^2\rho^2}
    H_{\hat x}(e_{k-1},t_{k-1})(y_k-y_{k-1})-\frac{1}{2\eps^2\rho^2}
    H_{\hat x}^2(e_{k-1},t_{k-1})\Delta t \r\}\right]\times\nonumber\\
    &&\nonumber\\
   && \left[\prod_{k=1}^N
\frac{\exp\l\{-\ds\frac{(y_k-y_{k-1})^2}{2\eps^2\rho^2\Delta t
}\r\}} {{\sqrt{2\pi\Delta t }}\,\varepsilon\rho}\right],\nonumber
 \eeq
where, by the Feynman-Kac formula \cite{Schulman}, \cite{Freidlin},
\cite{Kleinert}, \cite{unidirect}, \cite{ZakaiUstunel}, the first
product gives in the limit the function
 \beqq
 &&\varphi(e,t,\rho)=\\
 &&\\
 &&\lim_{N\to\infty} \underbrace{\int_{L}de_1\int_{L}de_2\cdots
\int_{L}de_{N-1}}_{N-1}\prod_{k=1}^N\left[\frac{1}{{\sqrt{2\pi\Delta
t }}\,\varepsilon\sigma}
\times\right.\\
&&\\
&&\left.\exp\l\{-\frac{[e_k-e_{k-1}-\Delta t M_{\hat x}
(e_{k-1},t_{k-1})\,]^2}{2\eps^2\sigma^2\Delta t }\r\}\times\right.\nonumber\\
   &&\nonumber\\
 &&\left. \exp\l\{\frac{1}{\eps^2\rho^2}
    H_{\hat x}(e_{k-1},t_{k-1})(y_k-y_{k-1})-\frac{1}{2\eps^2\rho^2}
    H_{\hat x}^2(e_{k-1},t_{k-1})\Delta t \r\}\right],
 \eeqq
which is the solution of the Zakai's equation in Stratonovich form
\cite{WZ}
 \beq
d_S\varphi(e,t,\rho)&=&\left\{-[\,M_{\hat
x}(e,t)\varphi(e,t)\,]_e+\frac{1}{2}[\,\eps^2\sigma^2
\varphi(e,t\,]_{ee}-\frac{\varphi(e,t)H_{\hat x}^2(e,t)}{2\eps^2\rho^2}\right\}\,dt +\nonumber\\
     &&\nonumber\\
    && \frac{\varphi(e,t)H_{\hat x}(e,t)}{\eps^2\rho^2}
    \,d_Sy(t),\label{ZakaiStr1}
 \eeq
with the boundary conditions
 \beq
 \varphi(e,t,\rho)=0\quad\mbox{for}\quad e\in\p L.\label{ZBC}
 \eeq
Therefore the joint density
 \beqq
&&p_N(e_N,t_N;\,y_1,y_2,\ldots,y_N)=\\
&&\\
&&
\Pr\{e_N(t_N)=e_N,\tau>t;y_N(t_1)=y_1,y_N(t_2)=y_2,\ldots,y_N(t_N)=y_N\}
 \eeqq
can be written  at $t=t_N, e_N=e$ as
 \beq
 p_N(e,t;\,y_1,y_2,\ldots,y_N)&=&\left[\varphi(e,t,\rho)+o(1)\right]
 \prod_{k=1}^N\frac{1}{\sqrt{2\pi\Delta t}\varepsilon\rho}
   \exp\l\{-\frac{(y_k-y_{k-1})^2}{2\eps^2\rho^2\Delta t}\r\},
   \label{jointN}
 \eeq
where $o(1)\to0$ as $N\to\infty$. Equivalently,
 \beq
   \varphi(e,t,\rho)=\frac{p_N(e,t;y_1,y_2,\ldots,y_N)}
   {\ds\prod_{k=1}^N\frac{1}{\sqrt{2\pi\Delta t}\varepsilon\rho}
   \exp\l\{-\frac{(y_k-y_{k-1})^2}{2\eps^2\rho^2\Delta t}\r\}}+o(1),
   \label{aposteriori}
   \eeq
which can be interpreted as follows: $\varphi(e,t,\rho)$ is the
joint conditional density of $e_N(t)$ and $\tau>t$, given the entire
trajectory $\{y_{N}(t_i)\}_{i=0}^N$, however, the probability
density of the trajectories $\{y_{N}(t_i)\}_{i=0}^N$,
 \beqq
 p_N^B(y_0^t)=\prod_{k=1}^N\left[\frac{\exp\l\{-\ds\frac{(y_k-y_{k-1})^2}
 {2\eps^2\rho^2\Delta t}\r\}}{\sqrt{2\pi\Delta
 t}\varepsilon\rho}\right],
   \eeqq
is Brownian, rather than the a priori density imposed by
(\ref{err_eq}), (\ref{err_meas}).

Now,
 \beqq
&&\Pr\{\tau>t_N,y_N(t_1)=y_1,y_N(t_2)=y_2,\ldots,y_N(t_N)=y_N\}=\\
&&\\
&&
\Pr\{\tau>t_N\,|\,y_N(t_1)=y_1,y_N(t_2)=y_2,\ldots,y_N(t_N)\}\times\\
&&\\
&& \Pr\{y_N(t_1)=y_1,y_N(t_2)=y_2,\ldots,y_N(t_N)=y_N\},
 \eeqq
which we abbreviate to
 \beq
 \Pr\{\tau>t,y_0^t\}=\Pr\{\tau>t\,|\,y_0^t\}p_N(y_0^t),\label{Abbreviated}
 \eeq
where the density
$p_N(y_0^t)=\Pr\{y_N(t_1)=y_1,y_N(t_2)=y_2,\ldots,y_N(t_N)=y_N\}$ is
defined by the system (\ref{dxN}), (\ref{dyN}), independently of
$\hat x(t)$.

We now use the abbreviated notation (\ref{Abbreviated}) to write
 \beq
 \Pr\{\tau>t\,|\,y_0^t\}&=&
 \frac{\Pr\{\tau>t,y_N(t_1)=y_1,y_N(t_2)=y_2,\ldots,y_N(t_N)=y_N\}}{p_N(y_0^t)}
 \nonumber\\
 &&\nonumber\\
&=&\int_{L}\frac{p_N(e,t;y_1,y_2,\ldots,y_N)}{p_N(y_0^t)}
\,de\nonumber\\
 &&\nonumber\\
&=&\frac{p_N^B(y_0^t)}{p_N(y_0^t)}\int_{L}\left\{\varphi(e,t,\rho)+o(1)\right\}\,de.
\label{Survival}
 \eeq
As $N\to\infty$, both sides of eq.(\ref{Survival}) converge to a
finite limit, which we write as
 \beq
 \Pr\{\tau>t\,|\,y_0^t\}=\alpha(t)\int_L\varphi(e,t)\,de,\label{SPHI}
 \eeq
where
$$\alpha(t)=\lim_{N\to\infty}\frac{p_N^B(y_0^t)}{p_N(y_0^t)},$$
is a function independent of $\hat x(t)$.

Next, we show that $E[\tau\,|\,y_0^\tau]$, as defined in
(\ref{EtauN}), (\ref{Etauy}), is given by
 \beq
E[\tau\,|\,y_0^\tau]=\int_0^\infty\Pr\{\tau>t\,|\,y_0^t\}\,dt.\label{EPtau}
 \eeq
Indeed, since $\Pr\{\tau>t\,|\,y_0^t\}\to0$ exponentially fast as
$t\to\infty$, we can write
 \beqq
\int_0^\infty\Pr\{\tau>t\,|\,y_0^t\}\,dt=\lim_{T\to\infty}\int_0^T
td\Pr\{\tau <t\,|\,y_0^t\}
 \eeqq
and
 \beqq
\int_0^T td\Pr\{\tau
<t\,|\,y_0^t\}=\lim_{N\to\infty}\sum_{i=1}^Ni\Delta t\Delta
\Pr\{\tau <i\Delta t\,|\,y_0^{i\Delta t}\},
 \eeqq
where
 \beqq
\Delta \Pr\{\tau <i\Delta t\,|\,y_0^{i\Delta t}\}= \Pr\{\tau
<i\Delta t\,|\,y_0^{i\Delta t}\}-\Pr\{\tau<(i-1)\Delta
t\,|\,y_0^{(i-1)\Delta t}\}.
 \eeqq
Now, we renumber the sampled
trajectories $e_{j,N}(t_i)$ in the numerator in (\ref{EtauN})
according to increasing $\tau_{i,N}$, so that in the new enumeration
$\tau_{i,N}=i\Delta t$. Then we group together the terms in the sum
that have the same $\tau_{i,N}$ and denote their sums $m_{i,N}$, so
that (\ref{EtauN}) becomes
 \beq
 E\left[\tau_N\wedge T\,|\,\left\{y_N(t_i),\,i=0,1,\ldots,\ds\left[\frac{\tau_{N}}
 {\Delta t}\right]\wedge N\right\}\right]&=&\frac{\ds\sum_{i=1}^Ni\Delta t\:m_{i,N}}
 {\ds\sum_{i=1}^N m_{i,N}}.
 \label{rearranged}
 \eeq
Finally, we identify
 \beqq
\Delta \Pr\{\tau <i\Delta t\,|\,y_0^{i\Delta
t}\}=\frac{m_{i,N}}{\ds\sum_{i=1}^{N}m_{i,N}}\left(1+o(1)\right)
 \eeqq
where $o(1)\to0$ as $N\to\infty$. Hence (\ref{EPtau}) follows.
Finally, we identify
 \beqq
\Delta \Pr\{\tau <i\Delta t\,|\,y_0^{i\Delta
t}\}=\frac{m_{i,N}}{\ds\sum_{i=1}^Nm_{i,N}}\left(1+o(1)\right)
 \eeqq
where $o(1)\to0$ as $N\to\infty$. Hence (\ref{EPtau}) follows.
\subsection{Asymptotic solution of Zakai's equation and the optimal filter}
For small $\varepsilon$ the solution of (\ref{ZakaiStr1}) with the
boundary conditions (\ref{ZBC}) is constructed by the method of
matched asymptotics \cite{Bender}, \cite{MS}, \cite{book}. The outer
solution is given by large deviations theory \cite{Hijab},
\cite{Freidlin}, \cite{Stroock}, \cite{DZ} as
 \beqq
    \varphi_{\mbox{\footnotesize { outer}}}(e,t)&=&\exp\left\{-\frac{\psi(e,t)}{\varepsilon^2}\right\},
 \eeqq
where
 \beq\label{psi}
    \psi(e,t)=\inf_{\ds e(\cdot)\in{\cal  C}_e^1([0,t])}\int_0^{t\wedge\tau}\left\{\left[\frac{\dot
    e(s)-M_{\hat x}(e(s),s)}{\sigma}\right]^2+\left[\frac{\dot
    y(s)-H_{\hat x}(e(s),s)}{\rho}\right]^2\right\}\,ds,
 \eeq
and
 \beqq
    {\cal C}_e^1([0,t])&=&\left\{e(\cdot)\in{\cal
    C}^1([0,t])\,:e(0)=e\right\}.
 \eeqq
We denote by $\tilde e(t)$ the minimizer of the integral on the
right hand side of eq.(\ref{psi}). The outer solution
$\varphi_{\mbox{\footnotesize { outer}}}(e,t)$ does not satisfy the
boundary conditions (\ref{ZBC}), so a boundary layer correction
$k(e,t,\varepsilon)$ is needed to obtain a uniform asymptotic
approximation,
 \beq
\varphi(e,t)\sim
\varphi_{\mbox{\footnotesize{uniform}}}(e,t)=\varphi_{\mbox{\footnotesize{outer}}}(e,t,\rho)k(e,t,\varepsilon)=
 \exp\left\{-\frac{\psi(e,t)}{\varepsilon^2}\right\}k(e,t,\varepsilon).\label{uniform}
 \eeq
The boundary layer function has to satisfy the boundary and matching
conditions
 \beq
  k(e,t,\varepsilon)=0\quad\mbox{for}\quad e\in\p L,\quad
\lim_{\varepsilon\to0}k(e,t,\varepsilon)=1\quad\mbox{for}\quad e\in
L,\label{match}
 \eeq
uniformly on compact subsets of the interior of $L$.

Since the survival probability is
 \beqq
 \Pr\left\{\tau>t\,|\,y_0^t\right\}
 =\int_{L}\alpha(t)\exp\left\{-\frac{\psi(e,t)}{\varepsilon^2}\right\}k(e,t,\varepsilon)\,de,
 \eeqq
the MTLL, according to  (\ref{EPtau}), is given by
 \beq\label{ETAU}
 E[\tau\,|\,y_0^\tau]=
 \int_0^\infty\int_{L}\alpha(t)\exp\left\{-\frac{\psi(e,t)}{\varepsilon^2}\right\}
 k(e,t,\varepsilon)\,de\,dt.
 \eeq

In view of (\ref{exhat}), the minimizer $\tilde e(t)$ of the
integral on the right hand side of (\ref{psi}) can be represented as
$\tilde e(t)=\tilde x(t)-\hat x(t)$, where $\tilde x(t)$ is the
minimizer of the integral
 \beq
    \Psi(x,t)=\inf_{\ds x(\cdot)\in{\cal
    C}_x^1([0,t])}\int_0^{t\wedge\tilde\tau}\left\{\left[\frac{\dot
    x(s)-m(x(s),s)}{\sigma}\right]^2+\left[\frac{\dot
    y(s)-h(x(s),s)}{\rho}\right]^2\right\}\,ds,\label{xtilde}
 \eeq
where $\tilde\tau=\inf\{t\,:\,\tilde x(t)-\hat x(t)\in\p L\}$
and
 \beqq
    {\cal C}_x^1([0,t])=\left\{x(\cdot)\in{\cal
    C}^1([0,t])\,:x(0)=x\right\}.
 \eeqq
Writing $\psi(e,t)=\Psi(x,t)$ and
$k(e,t,\varepsilon)=K(x,t,\varepsilon)$, we rewrite (\ref{ETAU}) as
 \beq
E[\tau\,|\,y_0^\tau]= \int_0^\infty\int_{L+\hat
x(t)}\alpha(t)\exp\left\{-\frac{\Psi(x,t)}{\varepsilon^2}\right\}
K(x,t,\varepsilon)\,dx\,dt.\label{ETAU2}
 \eeq
The integral in (\ref{ETAU2}) is evaluated for small $\varepsilon$
by the Laplace method, in which the integrand is approximated by a
Gaussian density with mean $\tilde x(t)$ and variance proportional
to $\varepsilon^2$. It is obviously maximized over the functions
$\hat x(t)$ by choosing $\hat x(t)$ so that the domain of
integration covers as much as possible of the area under the
Gaussian bell. If $L$ is an interval, then the choice $\hat
x(t)=\tilde x(t)$ is optimal. We conclude that for small noise, the
minimum noise energy filter $\tilde x(t)$ is asymptotically the
maximum MTLL filter.

\section{Discussion}

The main result of this paper is a proof that for small noise,
the minimum noise energy filter maximizes the mean time the
estimation error stays within a given region, e.g., maximizes the
mean time to lose lock in problems of phase tracking and
synchronization. The MNE filter is not finite-dimensional, however
finite discrete approximations, such as Viterbi-type algorithms
\cite{Unger}, \cite{ViterbiCoding}, can give arbitrary accuracy. The
practical aspects of finding the true MNE filter, or otherwise
adequate approximations for it, was partially dealt with in
\cite{Ehud} and still remains an interesting issue for further
studies.

Katzur {\em et. al.} \cite {KBS}, and subsequently Picard
\cite{Picard86}\cite{Picard91}, have shown that for nonlinear, but
monotone measurement functions, the MNE filter is to leading order
identical to the extended Kalman filter. However, for measurement
functions which are non-monotone, this is apparently not the case.
Ezri \cite{Doron} and Fischler \cite{Ehud} have considered the
problem of phase filtering and smoothing respectively, in which the
stochastic phase process $x(t)$ is measured in a low noise channel
by the vector function $\h(x) =
\left[\:\sin(x),\,\cos(x)\:\right]^T$. They show that there is a
huge gap between the MTLLs of the extended Kalman filter (smoother)
or particle filter, and the MNE filter (smoother), respectively.

The great advantage of the MNE filter in the case of phase
estimation is explained by the observation that finite-dimensional
approximations to the MAP or minimal MSEE filters (the EKF or the
finite dimensional filters of Katzur \cite{KBS}), do not capture
large deviations of the signal or of the measurements noise. They
are optimal only near local maxima of the a posteriori probability
density. The MNE filter, in contrast, is a \textit{global} MAP
estimator and can track large deviations. Thus, it is less
vulnerable to loss of lock phenomena, relative to the above
mentioned filters.\\

\noindent {\bf Acknowledgment:} The authors thank B.Z. Bobrovsky, O.
Zeitouni, D. Ezri, B. Nadler and A. Taflia for useful discussions.

\end{document}